\documentclass{article}

\usepackage{amsfonts}
\usepackage{amssymb}
\usepackage{amsfonts,amsmath}

\newtheorem{theorem}{Theorem}[section]
\newtheorem{definition}{Definition}[section]
\newtheorem{lemma}[theorem]{Lemma}
\newtheorem{proposition}[theorem]{Proposition}
\newtheorem{corollary}[theorem]{Corollary}

\newcommand{\eop }{ \hfill $\Box$ }

\begin{document}
\begin{center}
{\Large A characterization of Einstein manifolds\\}

\end{center}
\vspace{0.3cm}

\begin{center}
{\large  Sim\~ao Stelmastchuk}  \\

\textit{Departamento de  Matem\'{a}tica, FAFIUV \\ CX. P. 291 - CEP: 84600-000 - Uni\~ao da Vit\'oria - PR,\\
Brazil. e-mail: simnaos@gmail.com}
\end{center}

\vspace{0.3cm}

\begin{abstract}
In this work we wish characterize the Einstein manifolds $(M,g)$, however without the necessity of hypothesis of compactness over $M$ and unitary volume of $g$, which are well known in many works. Our result says that if all eingenvalues $\lambda$ of $r_{g}$, with respect to $g$, satisfy $\lambda \geq \frac{1}{n}s_{g}$, then $(M,g)$ is an Einstein manifold, where  $r_{g}$ and $s_{g}$ denote the Ricci and scalar curvatures, respectively.
\end{abstract}


\vspace{0.3cm} \noindent {\bf MSC2010 subject classification:}
53C25, 58J65, 60H30.

\section{Introduction}

Let $(M, g)$ be a Riemannian manifold. Let us denote the Ricci curvature by $r_{g}$ and the scalar curvature by
$s_{g}$. $M$ is called an Einstein manifold if, for every vector fields $X,Y$ on $M$, there exists a real constant such that
\[
r_{g}(X,Y) = \lambda g(X,Y).
\]
We call the metric $g$ of Einstein metric. We recall that scalar curvature $s_{g}$ is defined by $\mathrm{tr}\, r_{g}$.
A simple account gives that, for Einstein metrics, $s_{g}$ is the constant $\lambda n$, where $n$ is the dimension of
$M$. For a fuller treatment we refer the reader to \cite{besse}.

To characterize Einstein metrics we assume the fact that their scalar curvature is constant. However, we do not ask the
compactness property over Riemannian manifolds.

Our main result is proved using stochastic analysis on manifolds. The basic idea is, from hypothesis of scalar
curvature constat, show that the integral of the Ricci tensor $r_{g}- 1/n s_{g}g$ along any $g$-Brownian motion $B_{t}$
in $M$ is null, namely, $\int (r_{g}- 1/n s_{g}g)(dB,dB)=0$. From this we found conditions, see Proposition
\ref{proposition1}, to conclude that $r_{g}- 1/n s_{g}g$ is null, that is, $g$ is an Einstein metric. In summary, we
state our Theorem.\\

\emph{{\bf Theorem :} Let $(M,g)$ be a Riemannian manifold, $r_{g}$ its  Ricci curvature and $s_{g}$ its scalar
curvature. If all eigenvalues $\lambda$ of $r_{g}$, with respect to $g$, satisfy $\lambda \geq \frac{1}{n}s_{g}$, where
$n$ is the dimension of $M$, then $(M,g)$ is an Einstein manifold.}\\

This article is present in the following way: Section 2 contains a brief summary of Analysis Stochastic on Manifolds.
In Section 3 our main results are stated and proved.

\section{Stochastic tools}

In the following we always consider a complete probability space $(\Omega, \mathcal{F},\mathbb{P})$ endowed with a
filtration $(\mathcal{F}_{t})_{t \geq 0}$. We begin for introduce the three most important process for stochastic
analysis in manifolds. See for instance \cite{emery1} for a complete study about these process. From now on the term
smooth means of class $C^{\infty}$.

\begin{definition}
Let $M$ be a smooth manifold. A continuous $M$-valued process $X_{t}$ is called semimartingale if, for each smooth $f$
on $M$, the real-valued process $f\circ X_{t}$ is a semimartingale.
\end{definition}

Let $X_{t}$ be a semimartingale on $M$ and $b$ be a bilinear form on $M$. Let $(U, x_{1}, \ldots, x_{n})$ be a local
coordinate system on $M$. In this coordinate $b$ is written as $b_{ij} dx^{i} \otimes dx^{j}$, where $b_{ij}$ are
smooth function on $U$. The integral of $b$ along $X_{t}$ is defined, locally, by
\begin{equation}\label{eqtheorem2}
 \int b(dX,dX) = \int b_{ij}(X_{t}) d[X^{i},X^{j}]_{t},
\end{equation}
where $X^{i}_{t} = x^{i} \circ X_{t}$, $i=1, \ldots n$.

Using this definition has sense the following definition of martingales in smooth manifolds.

\begin{definition}
Let $M$ be a smooth manifold with a connection $\nabla$. A semimartingale $X_{t}$ in $M$ is called a martingale if, for
every smooth $f$,
\[
 f \circ X_{t} - f \circ X_{0} - \int \mathrm{Hess}f(dX,dX)
\]
is a real local martingale. Here, $\mathrm{Hess}$ denotes the Hessian operator associated to connection $\nabla$.
\end{definition}

In the sequel, we define Brownian motion in a smooth manifold.

\begin{definition}\label{brownian}
Let $(M,g)$ be a Riemannian manifold. Given $(\Omega,\mathcal{F},\mathbb{P},(\mathcal{F}_{t})_{t\geq 0})$, a $M$-valued
process $B_{t}$ is called a $g$-Brownian motion  in $(M,g)$ if $B_{t}$ is continuous and adapted and, for every smooth
f,
\[
 f \circ B_{t} - f \circ B_{0} - \frac{1}{2}\int \Delta_{g} f \circ B_{t} dt
\]
is a real local martingale, where $\Delta_{g}$ is the Laplace-Beltrami operator associated to metric $g$.
\end{definition}

Given a point $x$ in $(M, g)$, there always exists a $g$-Brownian motion $B_{t}$ in $M$, starting at $x$, defined on $[
0, \zeta [$ for some complete probability space $(\Omega, \mathcal{F},\mathbb{P},(\mathcal{F}_{t})_{t\geq 0})$ and some
stopping time $\zeta >0$.

%
%
%
%

The following Lemma, which demonstration is found in \cite[Lemma 5.20]{emery1}, is fundamental in the proof of our
Theorem.

\begin{lemma}\label{lemmafundmamental}
 If $B_{t}$ is a $g$-Brownian motion, then , for every bilinear form $b$,
\[
 \int b(dB,dB) = \int \mathrm{tr}\, b(B_{t}) dt.
\]
\end{lemma}


\section{Einstein manifolds}

We begin recalling the definition of Ricci tensor field for a metric and Einstein Manifold.

\begin{definition}
 The Ricci curvature $r_{g}$ of a Riemannian manifold $(M,g)$ is the 2-tensor
\[
 r_{g}(X,Y) = \mathrm{tr}\,(Z \rightarrow R(X,Z)Y),
\]
where $\mathrm{tr}$ denotes the trace of the linear map $Z \rightarrow R(X,Z)Y$.
\end{definition}

\begin{definition}
A Riemannian manifold $(M,g)$ is Einstein if there exists a real constant $\lambda$ such that
\[
 r_{g}(X,Y) = \lambda g(X,Y).
\]
We call the metric $g$ of Einstein metric.
\end{definition}

Before we prove our Theorem we show a general Proposition in Analysis Stochastic on Manifolds. This result is important
because in the proof of Theorem we have an step that is a particular case of this Proposition. In fact, this is the
final step to conclude the proof.

\begin{proposition}\label{proposition1}
Let $(M,g)$ be a Riemannian manifold. Let $b$ a symmetric bilinear form in $M$ such that their eigenvalues $\lambda_{i}
\geq 0$, $i=1, \ldots n$. If $\int b(dB,dB) = 0 $  for some $g$-Brownian motion $B_{t}$, then $b = 0$.
\end{proposition}
\begin{proof}
The proof is for a contrapositive argument. Let $b$ be a symmetric bilinear form on $M$ such that $b \geq 0$. Suppose
that $b>0$, that is, for $X,Y \in TM$, $b(X,Y) >0$. From Lemma \ref{lemmafundmamental}, for each $g$-Brownian motion
$B_{t}$, we obtain that
\[
\int b (dB, dB) = \int \mathrm{tr }b (B_{t}) dt = \sum_{i=1}^{n} \int \lambda_{i}(B_{t}) dt,
\]
We see that for each $\omega \in \Omega$, $\lambda_{i}(B_{t}(\omega))$ is continuous, for each $i=1, \ldots, n$. Since
$\lambda_{i}(B_{t})>0$ and $dt$ is the Lebesgue measure, it follows that $\int \lambda_{i}(B_{t}) dt > 0$. Thus we get
$\int b (dB, dB) >0$. \eop\\
\end{proof}



Now we prove our Theorem.

\begin{theorem}\label{teo1}
Let $(M,g)$ be a Riemannian manifold, $r_{g}$ its  Ricci curvature and $s_{g}$ its scalar curvature. If all eigenvalues $\lambda$ of $r_{g}$, with respect to $g$, satisfy $\lambda \geq \frac{1}{n}s_{g}$, where $n$ is the dimension of $M$, then $(M,g)$ is an Einstein manifold.
\end{theorem}
\begin{proof}
Let $(M,g)$ be a Riemannian manifold. Suppose that $s_{g}$ is constant, that is, there exists $ c \in \mathbb{R}$
such that $s_{g}(x) = c$ for all $x \in M$. If $n$ is dimension of $M$, then
\[
s_{g}(x) = \frac{\lambda}{n}n.
\]
As $n = \mathrm{tr}\, g$ we have
\[
 s_{g}(x) - \frac{c}{n} \mathrm{tr}\, g(x)= 0.
\]
Applying this equality about an arbitrary $g$-Brownian motion $B_{t}$ in $M$ we obtain
\[
 s_{g}(B_{t}) -\frac{c}{n}\mathrm{tr}\, g(B_{t}) = 0.
\]
We now integrate in $t$ each trajectory of $B_{t}$, that is,
\[
\int \mathrm{tr}(r_{g} - \frac{c}{n}g)(B_{t}) dt =\int \mathrm{tr}\,r_{g}(B_{t}) - \frac{c}{n}\mathrm{tr}\,
g(B_{t}) dt = 0.
\]
From Lemma \ref{lemmafundmamental} we conclude that
\begin{equation}\label{riccieq1}
\int (r_{g} - \frac{c}{n}g)(dB, dB) = 0.
\end{equation}
We now observe that eigenvalues of symmetric bilinear form $r_{g} - \frac{c}{n}g$ are $\lambda_{i}- \frac{c}{n}$, where $\lambda_{i}$, $i=1, \ldots, n$, are eigenvalues of $r_{g}$ with respect to $g$. By hypothesis, $\lambda_{i}- \frac{c}{n}\geq 0$. Using this fact in Proposition \ref{proposition1} we conclude that $r_{g} = \frac{c}{n}g$. Therefore, $g$ is an Einstein metric. \eop \\
\end{proof}

As a simple result from Theorem above we state an one about Ricci-Flat manifolds.

\begin{corollary}
Under hypothesis of Theorem \ref{teo1}, if $s_{g}=0$ then $(M,g)$ is Ricci-Flat manifold.
\end{corollary}

\end{document}